\documentstyle[12pt]{article}

\newtheorem{theorem}{Theorem}
\newtheorem{lemma}[theorem]{Lemma}
\newtheorem{prop}[theorem]{Proposition}
\newtheorem{corollary}[theorem]{Corollary}
\newtheorem{example}{Example}

\title{Semigroup cohomology and applications}
\author{B.\,V.\,Novikov (Kharkov, Ukraine)}
\date{}

\begin{document}
\maketitle

This article is a survey of the author's research. It consists of three sections
concerned three kinds of cohomologies of semigroups. Section 1 considers
`classic' cohomology as it was introduced by Eilenberg and  MacLane.
Here the attention is concentrated mainly on semigroups having cohomological
dimension 1. In Section 2 a generalization of the Eilenberg--MacLane cohomology
is introduced, the so-called 0-cohomology, which appears in applied topics
(projective representations of  semigroups, Brauer monoids). At last Section 3
is devoted to further generalizing: partial  cohomology defined and discussed
in it are used then for calculation of the classic cohomology for some
semigroups.

I am indebted to Prof. K.\,Roggenkamp and Prof. M.\,\c Stefanescu for the
support of my participation in the Workshop.

\section{EM-cohomology}

In this section we deal with the Eilenberg--MacLane cohomology of semigroups
\cite{c-e} (shortly EM-cohomology). Its definition is the same as for groups:
$$
H^n(S,A)={\rm Ext}^n_{{\bf Z}S}({\bf Z}, A).
$$
Here $S$ is a semigroup, $A$ a left $S$-module (i.\,e. a left module over the
integral semigroup ring ${\bf Z}S$), ${\bf Z}$ is considered as a trivial
$S$-module (i.\,e. $xa=a$ for all $x\in S$, $a\in {\bf Z}$); $H^n(S,A)$
is called {\it a $n^{th}$ cohomology group of $S$} (with the coefficient
module $A$).

Another definition (equivalent to preceding one) of $H^n(S,A)$ is following.
Denote by $C^n(S,A)$ the group of all maps
$f:\underbrace{S \times \ldots \times S}_{n\ {\rm times}} \to A$
(a group of $n$-cochains); a coboundary homomorphism
$\partial ^n:C^n(S,A)\to C^{n+1}(S,A)$ is given by the formula
\begin{eqnarray}\label{eq1}
\partial^nf(x_1,\ldots ,x_{n+1})=x_1f(x_2,\ldots ,x_{n+1}) \nonumber \\
+\sum_{i=1}^n{(-1)^if(x_1,\ldots ,x_ix_{i+1},\ldots ,x_{n+1})}
&+&(-1)^{n+1}f(x_1,\ldots ,x_n)
\end{eqnarray}
Then $\partial^{n}\partial^{n-1}=0$, i.\,e.
$$
{\rm Im}\partial^{n-1}=B^n(S,A)\ \mbox{(a group of coboundaries)}
$$
$$
\subseteq {\rm Ker}\partial^n=Z^n(S,A)\ \mbox{(a group of cocycles)}
$$
and cohomology groups are defined as $H^n(S,A)=Z^n(S,A)/B^n(S,A).$

It is worth noting two simple properties of the semigroup cohomology \cite{ber}
which are useful below.

(1) Adjoin to a given semigroup $S$ an extra element 1 and extend the
operation (multiplication) of $S$ to $S^1=S\cup \{1\}$ by
$$
\forall s\in S\quad s\cdot 1=1\cdot s=s, \qquad 1\cdot 1=1.
$$
Then $S^1$ becomes  a monoid (i.\,e. a semigroup with an identity element)
and is called {\it a semigroup with an adjoint identity}. Every $S$-module
turns naturally into a (unitary) $S^1$-module and we have
$$
H^n(S^1,A)\cong H^n(S,A), \qquad n\ge 0.
$$

(2) If $S$ possesses a zero element 0 then $H^n(S,A)=0$ for each $S$-module
$A$ and for  all $n\ge 1$. In particular, for any semigroup $S$ we can define
{\it the semigroup $S^0$ with an adjoint zero} (analogously to $S^1$); then
we have $H^n(S^0,A)=0$ ($n\ge 1$).

The semigroup EM-cohomology has not so wide applications as the cohomology
of groups. Nevertheless it is interesting for homologists at least as a model
for testing homological methods. The problem of describing semigroups
having cohomological dimension 1 is such an example. This problem has its
own story.

{\it Cohomological dimension} (c.\,d.) of a semigroup $S$ is the maximal
integer $n$ such that $H^n(S,A)\ne 0$ for some $S$-module $A$. There are many
reasons to study algebraic objects of c.d.\,1 (see, e.\,g., \cite{c-m}).

It is an easy exercise to prove that both a free group and a free
semigroup (or a free monoid) have c.d.\,1 \cite{c-e}. For groups the converse
is true --- this is the well-known Stallings--Swan theorem \cite{bro}.
So a group has cohomological dimension one if and only if it is free.

Now, what about  semigroups?

First, as we have mentioned above, every semigroup with 0 has c.d.\,$\le 1$.
This is the main reason why we have to confine ourselves to considering
cancellative semigroups only.

Second, a free group is not free as a semigroup. So even for cancellative
semigroups the Stallings--Swan theorem doesn't hold.

B.\,Mitchell \cite{mit}  has shown that a so-called partially free monoid
(the free product of a free group and a free monoid) has c.d.\,1. He has
supposed that if c.\,d.$S=1$ for  $S$ cancellative then $S$ is partially
free.

In \cite{nov3} I have built the first counter-example  for Mitchell
conjecture:
\begin{equation}\label{eq2}
S=\langle a,b,c,d\mid ab=cd\rangle
\end{equation}
and in \cite{nov4} I have formulated a `weakened Mitchell conjecture'. It
turned out true:

\begin{theorem} \label{1} {\rm  \cite{nov9}}
Every cancellative semigroup of c.d.\,$1$ can be embedded into a free group.
\end{theorem}

In the proof of this theorem the passage from the homological language to
the semigroup one is realized by the following lemma (which may be helpful
not only for semigroups).

Let $A$ be a left module over an arbitrary ring $R$,
$$
\ldots\mathop{\longrightarrow}\limits^{d_2}P_2
\mathop{\longrightarrow}\limits^{d_1}P_1
\mathop{\longrightarrow}\limits^{d_0}P_0
\mathop{\longrightarrow}\limits^{}0
$$
its projective resolution. Evidently, $d _{n}$ may be considered
as a $(n+1)$-dimensional cocycle with values in the $R$-module
${\rm Im\,} d_{n}$.

\begin{lemma}\label{2}
The cocycle $d_{n}\in Z^{n+1}(A, {\rm Im\,}d_{n})$
is a coboundary iff the projective dimension of $A$ is not
greater than $n$.
\end{lemma}

Applying Lemma \ref{2} to a bar-resolution of the $S$-module
${\bf Z}$ we obtain the next property of a cancellative semigroup $S$ with
c.d.\,1. Consider a graph with the elements of $S$ as vertices and with the
pairs $(a,b)\in S\times S$ such that $aS\cap bS\ne \emptyset$ as
edges. Then every circuit of this graph is triangulable.

This property allows to prove that $S$ can be embedded into a group (the latter
turns out free by the Stallings--Swan theorem). Note that
the converse assertion is certainly not true: there is a subsemigroup of a
free semigroup, which has c.d.\,1, while c.d. of its anti-isomorphic
is equal to 2 (see the examples below). A nice answer is only obtained in the
commutative case \cite{nov6}: the c.d. of a commutative cancellative
semigroup is equal to 1 if and only if this semigroup can be embedded into
${\bf Z}$.

\medskip

So a new problem arises: to describe subsemigroups of a free group having
c.d.\,1. This question seems rather difficult even  if we restrict ourselves
to subsemigroups of a free semigroup. The following results are taken out of
\cite{nov11}.

Let $S$ be a subsemigroup of a free semigroup $F$. Further development
of the proof of Theorem \ref{1} gives us the next assertion:

\begin{theorem}\label{3}
Let c.d.\,$S=1$ and $aS\cap S\ne\emptyset\ne Sa\cap S$ for some
$a\in F\setminus S$. There exists such $x\in aS\cap S$ that
$aS\cap S\subset xF^1$.
\end{theorem}

This theorem allows to build a lot of subsemigroups of $F$ having c.d.\,$>1$.

\begin{example}{\rm
Let $a, p,q,r$ are different elements of $F$ such that:

1) $\min(|a|, |p|,|q|,|r|)=|a|$ ($|a|$ denotes the length of  the word $a$),

2) $p$ and $q$ begin with different letters.

Then the subsemigroup $S=\langle p,q,r,ap,aq,ra\rangle$ has c.d.\,$>1$.}
\end{example}

From Theorem \ref{3} a solution of the proposed problem for left ideals
is obtained:

\begin{prop}\label{4}
A left ideal of a free semigroup has c.d.\,$1$ iff it is free.
\end{prop}

\begin{corollary}\label{5}
Every proper two-sided ideal of a free semigroup has c.d.\,$>1$.
\end{corollary}

For principal right ideals the situation is similar:

\begin{prop}\label{6}
A principal right ideal of a free semigroup has c.d.\,$1$ iff it is free.
\end{prop}

However for arbitrary right ideals  the analog of Proposition \ref{4} is not
true:
\begin{example}{\rm
Let $F=\langle a,b\rangle$ be a free semigroup. Then
$R=\{b,aba\}F^1$ is not free but c.\,d.$R=1$.}
\end{example}

By the way, from these results a counter-example to  another conjecture follows.
Yu.\,Drozd supposed that for any $S\in F$ either $S$ or the
antiisomorphic to $S$ has c.d.\,1. Consider the principal left ideal
$L=F^1aba$ in a free semigroup $F=\langle a,b\rangle$. It is not free since
its generators $aba, (ab)^2, (ab)^2a$, $(ab)^3$ obey the relation
$$
(ab)^2\cdot (ab)^2a = (ab)^3\cdot aba.
$$
By Proposition \ref{4} c.d.\,$L>1$. Of course its antiisomorphic $R=abaF^1$
is not free too and c.d.\,$R>1$ by Proposition \ref{6}. Hence the pair
$(L,R)$ gives a counter-example to the conjecture.

\section{0-cohomology}

In order to see how 0-cohomology appears let us  try to define
a projective representation of a semigroup.

Let $K$ be a field, $K^{\times}$ its multiplicative group, $n$ a positive
integer, $M(n,K)$ the semigroup of all $n\times n$ matrices over $K$.
Define an equivalence: for $A,B\in M(n,K)$
$$
A\sim B \Longleftrightarrow
\exists \lambda\in K^{\times}\ A= \lambda B.
$$
Then $\sim$ is a congruence on the semigroup $M(n,K)$ and we can consider
a factor semigroup $PM(n,K)=M(n,K)/\sim$, `the projective linear semigroup'.

Like for groups we call {\it a projective representation} of a given
semigroup  $S$ a homomorphism $\Gamma : S \to PM(n,K)$.

Fix an element in each $\sim$-class. Then $\Gamma$ induces a map
$\Gamma' : S \to M(n,K)$.  Now we can redefine a projective representation
of $S$: it is a map $\Gamma' : S \to M(n,K)$  such that

1) $\Gamma' (x)\Gamma' (y)=0 \Longleftrightarrow \Gamma' (xy)=0$,

2) $\Gamma' (x)\Gamma' (y)=\Gamma' (xy)\rho(x,y)$, where
$\rho :S\times S\to K^{\times}$ is a partial function defined on
the subset $\{(x,y)\mid \Gamma' (xy)\ne 0\}$.

Certainly $\rho$ yields the equation
$$
\rho(x,y)\rho(xy,z)=\rho(x,yz)\rho(y,z)
$$
for $\Gamma' (xyz)\ne 0$ and can be used as the corresponding
2-cocycle (like a factor system in Group Theory) excepting its partiality.
Therefore we must anew define suitable cohomology as follows.

\medskip

Let $S$ be an arbitrary semigroup with a zero. An Abelian group $A$ is called
{\it a 0-module} over $S$, if an action $(S\setminus \{0\})\times A\to A$
is defined which satisfies for all $s,t\in S\setminus \{0\},\ a,b\in A$ the
following
conditions:
$$
s(a+b)=sa+sb,
$$
$$
st\ne 0 \Longrightarrow s(ta)=(st)a.
$$

{\it A n-dimensional 0-cochain} is a partial $n$-place map from
$S$ to $A$ which is defined for all $n$-tuples \mbox{$(s_1,\ldots,s_n)$,}
such that $s_1\cdot \ldots \cdot s_n\ne 0$. A coboundary homomorphism is
given like for the usual cohomology by the formula (\ref{eq1}).
The equality $\partial^{n}\partial^{n-1}=0$ is valid too. We denote
$$
{\rm Im}\partial^{n-1}=B^n_0(S,A)\ \mbox{(a group of 0-coboundaries)}
$$
$$
\subseteq {\rm Ker}\partial^n=Z^n_0(S,A)\ \mbox{(a group of 0-cocycles)}
$$
and {\it 0-cohomology groups} are defined as
$H^n_0(S,A)=Z^n_0(S,A)/B^n_0(S,A).$

Note that for a semigroup $T^0=T\cup\{0\}$ with an adjointed zero
$$
H^n_0(T^0,A) \cong H^n(T,A),
$$
so the 0-cohomology may be considered as a generalization of the
Eilenberg--MacLane cohomology.

Properties of 0-cohomology are not considered here since they follow from the
properties of partial cohomologies (see Section 3).

Before returning to the  projective representations we need
a semigroup-theoretic construction, the so-called semilattice of groups
\cite{c-p}.

Let $\Lambda$ be a semilattice (i.\,e. a partially ordered set in which every
two elements $\lambda,\mu$ have the greatest lower bound $\lambda\mu$) and
let $\{G_{\lambda}\mid\lambda\in\Lambda\}$  be a family of disjoint groups.
For each pair $\lambda,\mu\in\Lambda$ such that $\lambda\ge\mu$, let
$\varphi^{\lambda}_{\mu}:G_{\lambda}\to G_{\mu}$ be a homomorphism. Suppose that

1) $\varphi^{\lambda}_{\lambda}$ is identical for every $\lambda\in\Lambda$,

2) $\varphi^{\lambda}_{\mu}\varphi^{\mu}_{\nu}=\varphi^{\lambda}_{\nu}$
for all ${\lambda}\ge{\mu}\ge{\nu}$.

Define a  multiplication on the set $T=\bigcup_{\lambda\in\Lambda}G_{\lambda}$
by the rule: if $x\in G_{\lambda}$, $y\in G_{\mu}$
$$
xy=(\varphi^{\lambda}_{\lambda\mu}x)(\varphi^{\mu}_{\lambda\mu}y).
$$
Then $T$ becomes a semigroup which is called {\it a semilattice of groups}.

Now return to the projective representations. Recall \cite{c-r} that if $S$
is a group then one defines an equivalence on the set of the factor systems
of $S$ (which corresponds to the equivalence of projective representations);
the factor set by this equivalence is a group ${\rm Sch}(S,K)$ which is called
{\it a Schur multiplicator} and describes (in some sense) all projective
representations of $S$ over  $K$. It is well-known that
${\rm Sch}(S,K) \cong H^2(S,K^{\times})$,  where $K^{\times}$ is considered
as a trivial $S$-module.

What will be for semigroups? In this situation ${\rm Sch}(S,K)$ is not
a group (more exactly, it becomes an inverse semigroup). Let $\Lambda$ be a
semilattice of all two-sided ideals of  $S$ (including $\emptyset$ and $S$)
with respect to the inclusion and the union as a greatest lower bound. Then
restriction of 0-cochains induces homomorphisms
$$
\varphi^{I}_{J}: H^n_0(S/I,K^{\times})\longrightarrow H^n_0(S/J,K^{\times})
$$
for  ideals $I\subseteq J$ and we have a semilattice of groups
$\bigcup_{I\in\Lambda}H^n_0(S/I,K^{\times})$ (here for $I=\emptyset$ we set
$H^n_0(S/\emptyset,K^{\times})=H^n(S,K^{\times})$). The next assertion was
proved in \cite{nov2}:

\begin{theorem}
For every semigroup $S$ and every field $K$
$$
{\rm Sch}(S,K) \cong \bigcup_{I\in\Lambda}H^2_0(S/I,K^{\times})
$$
\end{theorem}

Note that even if $0\not\in S$ we have to use 0-cohomology for describing of
${\rm Sch}(S,K)$.

\medskip

Another application of the 0-cohomology appears in connection with the Brauer
monoid. This notion was introduced by Haile, Larson and Sweedler
\cite{hai}, \cite{h-l-s} while they studied the so-called strongly primary
algebras (a generalization of central simple ones). I shall not give their
original definition which is rather complicated. But it turned out that the
Brauer monoid can be defined in terms of the 0-cohomology \cite{nov7}. To do
it one must introduce a new notion, a modification of a group.

By {\it a modification} $G(\ast)$ of a group $G$ we mean a semigroup on
the set $G^0=G\cup \{0\}$ with an operation $\ast$ such that $x\ast y$ is equal
either to $xy$ or to 0, while
$$
0\ast x=x\ast 0=0\ast 0=0
$$
and the identity of $G$ is the same for the semigroup $G(\ast)$.

In other words, to obtain a modification, one must erase the contents of
some  inputs in the multiplication table of $G$ and insert there zeros so that
the new operation would be associative.

Note some general properties of modifications. First, a modification of
$G$ satisfies the weak cancellation condition: from $x\ast z=y\ast z\ne 0$
it follows $x=y$ and analogously for left cancellation. Second, let $U$
be a subgroup of all invertible elements in $G(\ast)$. Then its complement
$I=G(\ast)\setminus U$ is a two-sided ideal. One can show that if $G$ is
finite, $I$ is nilpotent.

Let $S=G(\ast)$ and $T=G(\star)$ be modifications of $G$. It is  clear that
$S\cap T=G(\circ)$  is a modification too, where
$$
x\circ y\ne 0 \Longleftrightarrow x\ast y\ne 0\ne x\star y.
$$

We write $S\prec T$ if $x\star y=0$ implies $x\ast y=0$ for all $x,y\in G$.
Obviously, all modifications of $G$ constitute a semilattice $M(G)$:
the greatest lower bound in it is $S\cap T$.

Each $G$-module $A$ can be turned into a 0-module over a modification
$S$ in a natural way. Moreover, if $S\prec T$ then each 0-module over $T$
is  transformed into a 0-module over $S$. Therefore for $S\prec T$
a homomorphism  is defined
$$
\varphi^{T}_{S}: H^n_0(T,A)\longrightarrow H^n_0(S,A)
$$
and we obtain a semilattice of groups $\bigcup_{S\in M(G)}H^n_0(S,A)$.

In particular, let $L$ be a finite-dimensional normal extension of
a field $K$ with the Galois group $G$. Then $L^{\times}$ is a $G$-module.
We define {\it a (relative) Brauer monoid} as
$$
{\rm Br}(G,L)=\bigcup_{S\in M(G)}H^2_0(S,L^{\times})
$$
(the adjective `relative' will be omitted since in this article relative
Brauer monoids are only considered).

In the case when operation $\ast$ is defined in such a way that $x\ast y=xy$
for $x,y\ne 0$, we have
$$
H^2_0(S,L^{\times}) \cong H^2(G,L^{\times}),
$$
so the Brauer group is a subgroup of the Brauer monoid.

One can hope that the Brauer monoid will be useful. For example,
it is well-known that the Brauer group is trivial for any finite field
whereas the Brauer monoid is not trivial for each non-trivial field extension.

The Brauer monoid classifies strongly primary algebras over a field like
the Brauer group classifies division algebras.

The use of 0-cohomology allows us to split the study of the Brauer monoid
into two problems:

1) describing all modifications of a given finite group,

2) computing 0-cohomology of a modification.

Both of them seem rather difficult, especially the first. Its solution is
unknown even for cyclic groups. In \cite{nov8} some class of modifications
of simple cyclic groups is described. It  implies that the number
of modifications of the  group ${\bf Z}_p$ is $O(p^2)$. All multiplication
tables of the modifications $S_1,\ldots,S_{15}$
of ${\bf Z}_5$ (up to automorphisms of the group) are  shown in
Table\,\ref{fig1}; for ${\bf Z}_7$ their number equals 145.

\begin{figure}
\[
\begin{tabular}{|c||c|c|c|c|}
\hline
$S_1$ & $a^1$ & $a^2$ & $a^3$ & $a^4$ \\ \hline
$a^1$ & $0$ & $0$ & $0$ & $0$ \\ \hline
$a^2$ & $0$ & $a^4$ & $0$ & $a^1$ \\ \hline
$a^3$ & $0$ & $0$ & $a^1$ & $0$ \\ \hline
$a^4$ & $0$ & $a^1$ & $0$ & $0$ \\ \hline
\end{tabular}
\quad
\begin{tabular}{|c|c|c|c|c|}
\hline
$S_2$ & $a^1$ & $a^2$ & $a^3$ & $a^4$ \\ \hline
$a^1$ & $0$ & $0$ & $0$ & $0$ \\ \hline
$a^2$ & $0$ & $a^4$ & $0$ & $0$ \\ \hline
$a^3$ & $0$ & $0$ & $a^1$ & $0$ \\ \hline
$a^4$ & $0$ & $0$ & $0$ & $0$ \\ \hline
\end{tabular}
\quad
\begin{tabular}{|c|c|c|c|c|}
\hline
$S_3$ & $a^1$ & $a^2$ & $a^3$ & $a^4$ \\ \hline
$a^1$ & $0$ & $0$ & $0$ & $0$ \\ \hline
$a^2$ & $0$ & $0$ & $0$ & $a^1$ \\ \hline
$a^3$ & $0$ & $0$ & $a^1$ & $a^2$ \\ \hline
$a^4$ & $0$ & $a^1$ & $a^2$ & $a^3$ \\ \hline
\end{tabular}
\]
\[
\begin{tabular}{|c|c|c|c|c|}
\hline
$S_4$ & $a^1$ & $a^2$ & $a^3$ & $a^4$ \\ \hline
$a^1$ & $0$ & $0$ & $0$ & $0$ \\ \hline
$a^2$ & $0$ & $0$ & $0$ & $a^1$ \\ \hline
$a^3$ & $0$ & $0$ & $a^1$ & $0$ \\ \hline
$a^4$ & $0$ & $a^1$ & $0$ & $0$ \\ \hline
\end{tabular}
\quad
\begin{tabular}{|ccccc|}
\hline
\multicolumn{1}{|c|}{$S_5$} & \multicolumn{1}{c|}{$a^1$} &
\multicolumn{1}{c|}{$a^2$} & \multicolumn{1}{c|}{$a^3$} & $a^4$ \\ \hline
\multicolumn{1}{|c|}{$a^1$} & \multicolumn{1}{c|}{$0$} & \multicolumn{1}{c|}{%
$0$} & \multicolumn{1}{c|}{$0$} & $0$ \\ \hline
\multicolumn{1}{|c|}{$a^2$} & \multicolumn{1}{c|}{$0$} & \multicolumn{1}{c|}{%
$0$} & \multicolumn{1}{c|}{$0$} & $a^1$ \\ \hline
\multicolumn{1}{|c|}{$a^3$} & \multicolumn{1}{c|}{$0$} & \multicolumn{1}{c|}{%
$0$} & \multicolumn{1}{c|}{$a^1$} & $0$ \\ \hline
\multicolumn{1}{|c|}{$a^4$} & \multicolumn{1}{c|}{$0$} & \multicolumn{1}{c|}{%
$0$} & \multicolumn{1}{c|}{$0$} & $0$ \\ \hline
\end{tabular}
\quad
\begin{tabular}{|c|c|c|c|c|}
\hline
$S_6$ & $a^1$ & $a^2$ & $a^3$ & $a^4$ \\ \hline
$a^1$ & $0$ & $0$ & $0$ & $0$ \\ \hline
$a^2$ & $0$ & $0$ & $0$ & $0$ \\ \hline
$a^3$ & $0$ & $0$ & $a^1$ & $a^2$ \\ \hline
$a^4$ & $0$ & $0$ & $a^2$ & $0$ \\ \hline
\end{tabular}
\]
\[
\begin{tabular}{|c|c|c|c|c|}
\hline
$S_7$ & $a^1$ & $a^2$ & $a^3$ & $a^4$ \\ \hline
$a^1$ & $0$ & $0$ & $0$ & $0$ \\ \hline
$a^2$ & $0$ & $0$ & $0$ & $0$ \\ \hline
$a^3$ & $0$ & $0$ & $a^1$ & $0$ \\ \hline
$a^4$ & $0$ & $a^1$ & $0$ & $0$ \\ \hline
\end{tabular}
\quad
\begin{tabular}{|c|c|c|c|c|}
\hline
$S_8$ & $a^1$ & $a^2$ & $a^3$ & $a^4$ \\ \hline
$a^1$ & $0$ & $0$ & $0$ & $0$ \\ \hline
$a^2$ & $0$ & $0$ & $0$ & $0$ \\ \hline
$a^3$ & $0$ & $0$ & $a^1$ & $0$ \\ \hline
$a^4$ & $0$ & $0$ & $a^2$ & $0$ \\ \hline
\end{tabular}
\quad
\begin{tabular}{|c|c|c|c|c|}
\hline
$S_9$ & $a^1$ & $a^2$ & $a^3$ & $a^4$ \\ \hline
$a^1$ & $0$ & $0$ & $0$ & $0$ \\ \hline
$a^2$ & $0$ & $0$ & $0$ & $0$ \\ \hline
$a^3$ & $0$ & $0$ & $0$ & $a^2$ \\ \hline
$a^4$ & $0$ & $0$ & $a^2$ & $0$ \\ \hline
\end{tabular}
\]
\[
\begin{tabular}{|c|c|c|c|c|}
\hline
$S_{10}$ & $a^1$ & $a^2$ & $a^3$ & $a^4$ \\ \hline
$a^1$ & $0$ & $0$ & $0$ & $0$ \\ \hline
$a^2$ & $0$ & $0$ & $0$ & $0$ \\ \hline
$a^3$ & $0$ & $0$ & $0$ & $a^2$ \\ \hline
$a^4$ & $0$ & $0$ & $a^2$ & $a^3$ \\ \hline
\end{tabular}
\quad
\begin{tabular}{|c|c|c|c|c|}
\hline
$S_{11}$ & $a^1$ & $a^2$ & $a^3$ & $a^4$ \\ \hline
$a^1$ & $0$ & $0$ & $0$ & $0$ \\ \hline
$a^2$ & $0$ & $0$ & $0$ & $0$ \\ \hline
$a^3$ & $0$ & $0$ & $0$ & $0$ \\ \hline
$a^4$ & $0$ & $a^1$ & $0$ & $0$ \\ \hline
\end{tabular}
\quad
\begin{tabular}{|ccccc|}
\hline
\multicolumn{1}{|c|}{$S_{12}$} & \multicolumn{1}{c|}{$a^1$} &
\multicolumn{1}{c|}{$a^2$} & \multicolumn{1}{c|}{$a^3$} & $a^4$ \\ \hline
\multicolumn{1}{|c|}{$a^1$} & \multicolumn{1}{c|}{$0$} & \multicolumn{1}{c|}{%
$0$} & \multicolumn{1}{c|}{$0$} & $0$ \\ \hline
\multicolumn{1}{|c|}{$a^2$} & \multicolumn{1}{c|}{$0$} & \multicolumn{1}{c|}{%
$0$} & \multicolumn{1}{c|}{$0$} & $0$ \\ \hline
\multicolumn{1}{|c|}{$a^3$} & \multicolumn{1}{c|}{$0$} & \multicolumn{1}{c|}{%
$0$} & \multicolumn{1}{c|}{$0$} & $0$ \\ \hline
\multicolumn{1}{|c|}{$a^4$} & \multicolumn{1}{c|}{$0$} & \multicolumn{1}{c|}{%
$a^1$} & \multicolumn{1}{c|}{$0$} & $a^3$ \\ \hline
\end{tabular}
\]
\[
\begin{tabular}{|c|c|c|c|c|}
\hline
$S_{13}$ & $a^1$ & $a^2$ & $a^3$ & $a^4$ \\ \hline
$a^1$ & $0$ & $0$ & $0$ & $0$ \\ \hline
$a^2$ & $0$ & $0$ & $0$ & $0$ \\ \hline
$a^3$ & $0$ & $0$ & $0$ & $0$ \\ \hline
$a^4$ & $0$ & $0$ & $a^2$ & $0$ \\ \hline
\end{tabular}
\quad
\begin{tabular}{|c|c|c|c|c|}
\hline
$S_{14}$ & $a^1$ & $a^2$ & $a^3$ & $a^4$ \\ \hline
$a^1$ & $0$ & $0$ & $0$ & $0$ \\ \hline
$a^2$ & $0$ & $0$ & $0$ & $0$ \\ \hline
$a^3$ & $0$ & $0$ & $0$ & $0$ \\ \hline
$a^4$ & $0$ & $0$ & $0$ & $0$ \\ \hline
\end{tabular}
\quad
\begin{tabular}{|c|c|c|c|c|}
\hline
$S_{15}$ & $a^1$ & $a^2$ & $a^3$ & $a^4$ \\ \hline
$a^1$ & $0$ & $0$ & $0$ & $0$ \\ \hline
$a^2$ & $0$ & $0$ & $0$ & $0$ \\ \hline
$a^3$ & $0$ & $0$ & $0$ & $0$ \\ \hline
$a^4$ & $0$ & $0$ & $0$ & $a^3$ \\ \hline
\end{tabular}
\]
\caption{The modifications of ${\bf Z}_5$}\label{fig1}
\end{figure}

As to the second problem, the initial step in solving it may consist in
eliminating the influence of invertible elements of modifications on
the structure of the Brauer monoid. Some results in this direction were
obtained in \cite{ki-n}, \cite{nov7}.

As above let $G$ be the Galois group of a finite-dimensional extension $L/K$,
$S=G(\ast)$ its modification, $U$ the subgroup of invertible elements of $S$,
$P$ the subfield of all $U$-fixed elements: $P=\{a\in L\mid Ua=a\}.$
The inclusion $U\hookrightarrow S$ induces a homomorphism
$$
\psi: H_0^2(S,L^{\times})\longrightarrow H^2(U,L^{\times})
$$
We shall study this homomorphism in the situation when $U$ is a normal subgroup
of $S$ (i.\,e. $x\ast U=U\ast x$ for all $x\in S$). Then the factor semigroup
$S/U$ is well-defined.

Further, if $U\triangleleft S$ then $P^{\times}$ is a 0-module over $S/U$.
The inclusion $P^{\times}\hookrightarrow L^{\times}$ and the epimorphism
$S\to S/U$ induce a homomorphism
$$
\chi: H_0^2(S/U,P^{\times})\longrightarrow H_0^2(S,L^{\times})
$$

\begin{theorem}\label{7}
Let $U\triangleleft S$. Then the sequence
$$
0\longrightarrow H^2_0(S/U, P^{\times})
\mathop{\longrightarrow}\limits^{\chi}
H_0^2(S,L^{\times})\mathop{\longrightarrow}\limits^{\psi} H^2(U,L^{\times})
$$
is exact.
\end{theorem}

\begin{corollary} \label{8}
If the field $L$ is finite then
$$
H_0^2(S,L^{\times})\cong H^2_0(S/U, P^{\times})
$$
\end{corollary}

Therefore, for finite fields the problem is reduced to the computation of
0-cohomology of a nilpotent 0-cancellative semigroup $(S/U)\setminus \{1\}$.
I believe that such an algorithm can be built.

At last note that W.\,Clark \cite{cla} used 0-cohomology (however with
trivial 0-modules only) for investigation of some matrix algebras.

\section{Partial cohomologies}

0-Cohomology has one more application: for calculating of EM-cohomology.
However from this point of view it is worth once more to generalize our
construction.

One can ask: what would be if we considered partial maps as cochains, starting
from an arbitrary subset $W\subseteq S$, not necessary from $S\setminus \{0\}$?
It was shown in \cite{nov4} that this question is reduced to the following
particular case.

Let a semigroup $S$ be generated by a subset $W$ with defining relations of the
form $xy=z$ for some $x,y,z\in W$. Such a $W$ will be called {\it a root} of
$S$. We denote by $W_n$ a set of all $n$-tuples $(x_1,\ldots,x_n)$ such
that $x_ix_{i+1}\cdot\ldots\cdot x_j \in W$ for all $1\le i\le j\le n$.
Every map from $W_n$ to a $S$-module $A$ is called {\it a partial
$n$-dimensional cochain} of $W$ or {\it a $W$-cochain}
with values in $A$. $n$-Dimensional $W$-cochains form an Abelian group
$C^n(S,W,A)$ for $n>0$. We set $C^0(S,W,A)=A$, and if $W_n=\emptyset$
then $C^n(S,W,A)=0$. The coboundary homomorphism is given by the same formula
(\ref{eq1}); the corresponding {\it partial cohomology groups}
(or {\it $W$-cohomology groups}) are denoted by $H^n(S,W,A)$.

It is clear that we obtain EM-cohomology if $W=S$. Reducing
0-cohomology to a partial one looks more complicated: if $S$ is a semigroup
with 0, $W=S\setminus \{0\}$, then we generate a new semigroup
$T=\langle W\rangle$ with  the operation $\ast$ and defining relations of
the form $u\ast v=w$, where $u,v,w\in W$ and $uv=w$ in $S$. Then $W$ is
a root in $T$ and
$$
H^n_0(S,A)\cong H^n(T,W,A).
$$

Having a presentation of a semigroup $S$ one can easily build some of its roots.

\begin{example}\label{ex}{\rm
Let $S=\langle a,b,c,d\mid ab=cd\rangle$ (see (\ref{eq2})). Then
$$
W=\{a,b,c,d,x=ab\}, \quad W_2=\{(a,b),(c,d)\}, \quad W_3=\emptyset
$$
and
$$
S=\langle a,b,c,d,x\mid ab=x,cd=x\rangle
$$}
\end{example}

How are $H^n(S,W,A)$ and $H^n(S,A)$ connected?

The embedding $W\hookrightarrow S$ induces homomorphisms
$$
\theta ^n_W:H^n(S,A)\longrightarrow H^n(S,W,A),
$$
\begin{prop}\label{9}{\rm \cite{nov5}}
If $W$ is a root of a semigroup $S$ then $\theta ^n_W$ is an isomorphism for
$n\le 1$ and a monomorphism for $n=2$.
\end{prop}

Generally speaking, $\theta^2_W$ can be non-surjective (by the way it means
that partial cohomology ought not be a derived functor in the category of
$S$-modules).

Proposition \ref{9} enables us to use partial cohomology for
calculating 1-dimensional EM-cohomology of semigroups (and getting some
information about 2-dimensional one) in the case when one succeeds to
find a `good' root in a given semigroup. For instance, consider
the semigroup $S$ from Example \ref{ex}. Define for each $f\in Z^2(S,W,A)$
the 1-dimensional $W$-cochain $h$ by
\begin{equation}
h(s)=
\left \{ \begin {array}{ll}
f(a,b), & {\rm if}\ s=a, \\
f(c,d), & {\rm if}\ s=c, \\
0,      & {\rm otherwise.}
                  \end {array} \right.
\end{equation}
Then $f=\partial\,h$, so $H^2(S,W,A)=0$. Since $\theta^2_W$ is injective,
$H^2(S,A)=0$ too.

To study $\theta ^n$ for $n>1$ we need some new definitions.

Let $S$ be a semigroup, $W$ be its root. A decomposition $x=x_1\ldots x_k $
($x_i\in W$) of an element $x\in S\setminus W$ is called {\it reduced}
if $x_ix_{i+1}\ldots x_j\not \in W$ for each $i,j,\ 1\le i< j\le k$. We
mean that a reduced decomposition of an element $x\in W$ is its
decomposition into the product of one multiplier. A root $W$ is said to be
{\it canonic} if each element $x\in S$ has the unique reduced decomposition.

For example, the set of all element of $S$ is a canonic root.

A root $W$ is called a {\it $J$-root} if $xy=x,\ yz=z$ implies $xz\in W$ 
for all $x,y,z\in W$.

\begin{theorem}\label{10}{\rm \cite{nov5,nov10}}
If $W$ is a canonic $J$-root of $S$, then $\theta^n_W$ are isomorphisms
for all $n\ge 0$.
\end{theorem}

As above, we can use Theorem \ref{10} for calculating EM-cohomology
in higher dimensions. For example, if $S=T*U$ is the free product of
semigroups $T$ and $U$, then  $W=T\cup U$ is its canonic $J$-root and
$W_n=T_n\cup U_n$. Thus, we get
$$
H^n(S,A)\cong H^n(S,W,A)\cong H^n(T,A)\oplus H^n(U,A)
$$
for every $S$-module $A$. Below we consider less trivial examples.

\begin{example}\label{ex1}{\rm
Let $S=\langle a,b_1,b_2,\ldots\mid aP=Q\rangle$ be such a semigroup that
the words $P$ and $Q$ do not contain the letter $a$. Denote
$W=F\cup\{a\}$, where $F=\langle b_1,b_2,\ldots\rangle$ is a subsemigroup
of $S$. Then $W$ is a canonic $J$-root. The fact that $F$ is free facilitates
the calculation of $W$-cohomology of $S$; so we obtain
$H^2(S,A)=0$ for every $X$-module $A$ (that is c.d.\,$S= 1$).}
\end{example}

\begin{example}{\rm
The semigroup
$$
S^{\rm op}=\langle a,b_1,b_2,\ldots\mid Pa=Q\rangle,
$$
is antiisomorphic to $S$ (see Example \ref{ex1}). Like for $S$, the subset
$W=F\cup\{a\}$ is a canonic $J$-root. However c.d.\,$S^{\rm op}= 2$.
Besides,  $H^2(S^{\rm op},A)\cong A/B$, where
$$
B=PA+\sum_i\left(\frac{\partial P}{\partial b_i}-
\frac{\partial Q}{\partial b_i}\right)A;
$$
here $\frac{\partial\ }{\partial b}$ is an analog of the Fox' derivative
\cite{c-f} adapted to semigroups in \cite{nov5}.}
\end{example}

Consider one more pair of antiisomorphic semigroups.

\begin{example}{\rm
Let $U$ be an arbitrary semigroup,
$$
T=\langle U,p\mid Up=p\rangle \quad (p\not\in U)
$$
This notation means that $T$ is generated by its subsemigroup $U$ and
by an element $p\not\in U$ and is defined by relations of the form
$$
u\cdot v =uv,\quad u\cdot p =p \quad (u,v\in U)
$$
The subset $W=U\cup\{p\}$ turns out  a canonic $J$-root.
We get c.d.\,$T=1$ and $H^1(T,A)\cong A/(p-1)A$ for every
$T$-module $A$.}
\end{example}

\begin{example}{\rm
Now consider the semigroup
$$
T^{\rm op}=\langle U,p\mid pU=p\rangle,
$$
antiisomorphic to $T$. Its EM-cohomology is much more complicated:}
\end{example}

\begin{prop}\label{11}
Let $A$ be a $T^{\rm op}$-module, $A_1$ be its additive group considered as
a trivial $T^{\rm op}$-module. The homomorphisms
$\psi^n :H^n(T^{\rm op},A)\to H^n(U,A)$ induced by the embedding
$U\hookrightarrow T^{\rm op}$ are inserted into the long exact sequence
\begin{eqnarray}\label{eq3}
0\to H^0(T^{\rm op},A)\mathop{\to}\limits^{\psi^0}
H^0(U,A)\to H^0(U,A_1)\to H^1(T^{\rm op},A)\mathop{\to}\limits^{\psi^1} ...
\nonumber\\
...\to H^n(T^{\rm op},A)\mathop{\to}\limits^{\psi^n}
H^n(U,A)\to H^n(U,A_1)\to ...
\end{eqnarray}
\end{prop}

By the last two examples we can build a semigroup $T$ such that c.d.\,$T=1$
and c.d.\,$T^{\rm op}= \infty$. To do it take the additive group of
the  ring ${\bf Z}_9$ as $A$ and its multiplicative group as $U$.
The action of $U$ on $A$ coincides with the multiplication in ${\bf Z}_9$. Then
$H^n(T^{\rm op},A)\cong {\bf Z}_3$ for $n>1$.

\medskip

It is worth to add that the notion of a canonic root can be applied to
algorithmic  problems. For instance, with its help a new family of semigroups
with solvable word problem was obtained \cite{ka-n}.

\bigskip

\bigskip

B.V.Novikov, Saltovskoye shosse 258, apt.20, Kharkov, 68178, Ukraine

e-mail: boris.v.novikov@univer.kharkov.ua


\begin{thebibliography}{99}

\bibitem{ber}
N.\,Bernstein.
        \newblock{\it On the cohomology of semigroups.}
        \newblock{Dissert. Abstrs., {\bf 25}(1965), N1, 6644-6645.}
\bibitem{bro}
K.\,S.\,Brown.
        \newblock{\it Cohomology of Groups.}
        \newblock{Springer-Verlag, 1982.}
\bibitem{c-e}
H.\,Cartan and S.\,Eilenberg.
        \newblock{\it Homological Algebra.}
        \newblock{Princeton, 1956.}
\bibitem{c-m}
C.\,C.\,Cheng and B.\,Mitchell.
        \newblock{\it  DCC posets of cohomological dimension one.}
        \newblock{J. Pure and Appl. Algebra, {\bf 13}(1978), N2, 125-137.}
\bibitem{cla}
W.\,E.\,Clark.
        \newblock{\it Cohomology of semigroups via topology with an
        application to semigroup algebras.}
        \newblock{Commun. Algebra, {\bf 4}(1976), 979-997.}
\bibitem{c-p}
A.\,H.\,Clifford and G.\,B.\,Preston.
        \newblock{\it The Algebraic Theory of Semigroups.}
        \newblock{vol.1, 2, AMS Math. Surveys, 1964, 1967.}
\bibitem{c-f}
R.\,H.\,Crowell and R.\,H.\,Fox.
        \newblock{\it Introduction to Knot Theory.}
        \newblock{Ginn \& Co, 1967.}
\bibitem{c-r}
        C.\,W.\,Curtis and I.\,Reiner.
        \newblock{\it Representation Theory of Finite Groups
        and Associative Algebras.}
        \newblock{Willey \& Sons, N.-Y., L., 1967.}
\bibitem{hai}
D.\,E.\,Haile.
        \newblock{\it The Brauer monoid of a field. }
        \newblock{J. Algebra, {\bf 81}(1983), N2, 521-539.}
\bibitem{h-l-s}
D.\,E.\,Haile, R.\,G.\,Larson and M.\,E.\,Sweedler.
        \newblock{\it A new invariant for ${\bf C}$ over ${\bf R}$: almost
        invertible cohomology theory and the classification of idempotent
        cohomology classes and algebras by partially ordered sets
        with Galois group action.}
        \newblock{Amer. J. Math., {\bf 105}(1983), N3, 689-814.}
\bibitem{ka-n}
O.\,S.\,Kashcheeva and B.\,V.\,Novikov.
        \newblock{\it Canonic subsets in semigroups.}
        \newblock{Filomat (Yugoslavia), {\bf 12}(1998), N1, 21-27.}
\bibitem{ki-n}
V.\,V. Kirichenko and B.\,V.\,Novikov.
       \newblock{\it On the Brauer monoid for finite fields.}
       \newblock{In: Proc. of 5th Internat. Conf. on Finite Fields and Appl.
       (to appear).}
\bibitem{mit}
B.\,Mitchell.
	\newblock{\it On the dimension of objects and
        categories. {\rm I.} Monoids.}
        \newblock{J. Algebra, {\bf 9}(1968), N3, 314--340.}
\bibitem{nov1}
B.\,V.\,Novikov.
        \newblock{\it On 0-cohomology of semigroups.}
        \newblock{In: ``Theor. and Appl. Quest. of Diff. Equat. and Algebra".
        Kiev, Naukova dumka, 1978, 185-188 (in Russian).}
\bibitem{nov2}
B.\,V.\,Novikov.
        \newblock{\it On projective representations of semigroups.}
        \newblock{Doklady AN USSR, 1979, N6, 474-478 (in Russian).}
\bibitem{nov3}
B.\,V.\,Novikov.
        \newblock{\it An counter-example to a Mitchell conjecture.}
        \newblock{Trudy Tbiliss. matem. inst. AN GSSR, {\bf 70}(1982),
        52-55 (in Russian).}
\bibitem{nov4}
B.\,V.\,Novikov.
        \newblock{\it On partial cohomologies of semigroups.}
        \newblock{Semigroup Forum, {\bf 28}(1984), N1-3, 355-364.}
\bibitem{nov5}
B.\,V.\,Novikov.
        \newblock{\it Partial cohomologies and their applications.}
        \newblock{Izv. vuzov. Matem., 1988, N11, 25-32 (in Russian).
        Translation in Soviet Math. (Iz. VUZ) v.32 (1988), N11, 38-48.}
\bibitem{nov6}
B.\,V.\,Novikov.
        \newblock{\it Commutative cancellative semigroups of dimension 1.}
        \newblock{Matem. Zametki, {\bf 48}(1990), N1, 148-149 (in Russian).}
\bibitem{nov7}
B.\,V.\,Novikov.
        \newblock{\it On the Brauer monoid.}
        \newblock{Matem. Zametki, {\bf 57}(1995), N4, 633-636 (in Russian).
        Translation in Math. Notes 57 (1995), no. 3-4, 440-442.}
\bibitem{nov8}
B.\,V.\,Novikov.
        \newblock{\it On modification of the Galois group.}
        \newblock{Filomat (Yugoslavia), {\bf 9}(1995), N3, 867-872.}
\bibitem{nov9}
B.\,V.\,Novikov.
       \newblock{\it Semigroups of cohomological dimension 1.}
       \newblock{J. Algebra, {\bf 204}(1998), 386-393.}
\bibitem{nov10}
B.\,V.\,Novikov.
       \newblock{\it Partial cohomologies and canonic roots in semigroups.}
       \newblock{Matem. studii (Ukraine), {\bf 12}(1999), N1, 7-14 (in Russian).}
\bibitem{nov11}
B.\,V.\,Novikov.
       \newblock{\it On cohomological dimension of ideals of free semigroups.}
       \newblock{In: Colloq. on Semigroups. July 17-21, 2000. Szeged.
       Abstracts, 19.}

\end{thebibliography}
\end{document}